\documentclass{amsart}
\usepackage{amssymb, array, verbatim, amscd}
\usepackage{amsmath, amsfonts,amsthm}

\theoremstyle{plain}
\newtheorem{theorem}{Theorem}[section]

\theoremstyle{definition}

\newtheorem{definition}[theorem]{Definition}
\newtheorem{example}[theorem]{Example}

\newtheorem{remark}[theorem]{Remark}

\newcounter{abc}

\newcounter{ABC}

\newcommand{\Z}{{\mathbb Z}}

\newcommand{\R}{{\mathbb R}}

\newcommand{\F}{{\mathcal F}}

\begin{document}
\title[Genetic Sequential  Dynamical Systems ] {Genetic Sequential
Dynamical Systems }

\author{M. A. Avi\~no}
\address{Department of Physics-Mathematics\\
University of  Puerto Rico\\
Cayey, PR  00736 } \email{m\_avino@cayey1.upr.clu.edu}
 \author{H. Ortiz}
 \address{Programmer-Archaeologist High Performance Computing
 facility\\
University of Puerto Rico\\ http://www.hpcf.upr.edu/~humberto/ }
 \email{humberto@hpcf.upr.edu }
\author{O. Moreno}
 \address{Department of Mathematics, and Computer Sciences\\
University of Puerto Rico }
 \email{moreno@uprr.pr }



\begin{abstract}
The whole complex process to obtain a protein encoded by a gene is
difficult to include in a mathematical model. There are many
models for describing different aspects of a genetic network.
Finding a better model is one of the most important and
interesting questions in computational biology. Sequential
dynamical systems have been developed for a theory of computer
simulation, and in this paper, a genetic sequential dynamical
system is introduced. A gene is considered to be a function which
can take a finite number of values. We prove that a genetic
sequential dynamical system is a mathematical good description for
a finite state linear model introduced by Brazma \cite{B,B2}.
\end{abstract}

\maketitle

  \section{Introduction}
    One important  and interesting question in biology is how
     genes are regulated. The most important models for gene
     regulation networks are boolean models and differential
     equation based models.
  Boolean models \cite{K,K1,A} describe the activity of genes using an element
  of $Z_2=\{0,1\}^n$, that is an $n$ dimensional  vector with entries
  in $\{0,1\}$. Each entry $x_i$ means the activation of the  gene
$i$. In the Boolean model   we have  the vector
  space ${\Z_2}^n$ with an attached  function $f:{\Z_2}^n\rightarrow
  {\Z_2}^n$. The iteration of $f$ means the time passed. The
  properties of the digraph associated with these iterations
   are the characteristics of the network.  There are different ways to
generalize this model:  using more than two possibilities for each
gene, and second using several functions for each gene, (PBN).
Recently, a new mathematical model \emph{Probabilistic Boolean
Networks} (PBN) was introduced by Ilya Shmulevich, \cite{SDZ}.
This model introduced the probabilistic behavior in Boolean
networks and has been  used to predict the steady states of
genetic networks in cancer cells, \cite{I4}. For other
mathematical models see\cite{J}.

If we assume that a gene has more than two levels of possibilities
which are determined  by the environment and the concentration of
a particular substance in the network, then  the activity of  a
gene $i$ is taken from the set of natural numbers $\{0,1,\ldots
,m-1\}$, where $m>1$, \cite{T,T1}. One of the problem to study
genetic networks with more than two possibilities for each gene is
to find a good way to describe the functions in the net. There are
polynomial representations for this functions over a finite field
and one of the most important results using the techniques
  of Computer Algebra is that we can find all the polynomial solutions.
  This ideas appeared  for a first time in the Seminar of Reinhard
   Laubenbacher in Virginia Bioinformatic Tech \cite{LS, LS2}. In
   \cite{A1}, these functions are called partially defined
   functions and study over a finite field of $p^n$ elements.
   In \cite {G, G2}, it  is proved that there exist
   polynomials solutions.

 The theory of sequential dynamical systems (SDS) was first  introduced in \cite{BR, BRM1, BRM2}
  as a mathematical abstraction of a simulated system by a computer. In \cite {LP1, LP2}, Laubenbacher and
  Pareigis introduced a categorical framework for the study of SDS. In this paper, we describe a particular SDS for genetic
 networks. In addition, we present  a mathematical background to use in the
study of the Finite State Linear Model introduced by Brazma and
Schlitt,\cite{B2}. Using partially defined function and the
polynomial solutions we present a mathematical particular example
over a finite field of three elements.

 This paper is organized as follow, in section $2$, we describe  the  ideas of the
Finite State Linear Model  for Gene Regulation Networks introduced
by Alvis Brazma in \cite{B,B2}, and introduce a notation slightly
different from the one used in \cite {B,B2}. In section $3$, we
compare the  definitions of SDS (\cite{LP2}), and the definition
of parallel dynamical systems. In section $4$, we define two
different models: the first  describes the mathematical aspects of
the Brazma Model and the second (genetic sequential dynamical
system) generalizes the Brazma model.
      \section{Finite state linear model}
      Gene expression is a two-step process: first, a
      single stranded  messenger RNA (mRNA) is copied
      (transcribed) from the strand of a duplex DNA
       molecule that encodes genetic information. In the second
       step, the mRNA moves to the cytoplasm, is complexed to
       ribosomes, and its genetic information is translated into
       the amino acid sequence of a polypeptide.

The model (FSLM) in \cite{B2}  considers the following definitions
of a gene ( Section 3.3,   \cite{BPS}) and gene regulation
networks (Section 4.3, \cite{BPS}).
\begin{definition}\label{gen}
\emph{A \textbf{gene} is a continuous stretch of a genomic DNA
molecule, from which a complex molecular machinery can read
information (encoded as a string of $A$, $T$, $G$, and $C$) and
make a particular type of a protein or a few different proteins.}
\end{definition}

\emph{Transcription factors control a gene expression by binding
the gene's promoter and either activating (switching on) the
gene's transcription, or repressing it (switching off).
Transcription factors are gene products themselves, and therefore
in turn can be controlled by other transcription factors.
Transcription factors can control many genes, and some (probably
most) genes are controlled by combinations of transcription
factors. Feedback loops are possible. Therefore we can talk about
\textbf{gene regulation networks}. Understanding, describing and
modelling such gene regulation networks are one of the most
challenging problems in functional genomics.}
\bigskip

   Now, we introduce the notation that we use in this paper.
   The  network has  $n$ genes
      $G_1, \ldots , G_n$. The binding sites are stages of the
      processes of transcription of a particular gene $G_j$,
      denoted by $B_{j1},\ldots , B_{jm_j}$.
      Each binding site $B_{jk}$ is determined by the concentration
      $c_{j}(t)$ at time $t$,  of a particular  substance $i_{j}$
       associated with or generated by  a gene $G_j$.
    There are two  constants for each state of the binding site $B_{jk}$, $a_{jk}$
        and $d_{jk}$, called  association and dissociation constants
         respectively. That is, taking the real number
         $c_{jk}(t)$ and depending on its relation to $a_{jk}$
        and $d_{jk}$, we give a state for the binding sites
        $B_{jk}$.
        Each binding site  $B_{jk}$ can  take a finite number of
        possibilities. For  understanding the problem,  $B_{jk}$
         is taking as a finite subset of the set of the integers $\Z$.
        Let denote  by $b_{jk}\in \R$ a state of
        $B_{jk}$. In FSLM, they called $B_{jk}$ a multistate binding site,
        and the  vector $b_j=(b_{j1},\ldots ,b_{jm})$ a binding site
vector.
        The set of all possible vectors $b_j$ is the environment of
        the gene $G_j$, that is
        \[B_j=B_{j1}\times \cdots \times B_{jm_j}\subset \Z ^{m_j}.\]

For each gene, there is a function $F_j$ (the control function).
Its inputs are  the
 states of the binding sites, $B_j$.
The function $F_j$ takes one of the values of the binding sites,
but
 the output is the production of a substance at a given rate which
can act again over all the binding sites.
 Therefore we consider the control functions $F_j:B_j
 \rightarrow B_j$, whose  output, a vector of $B_j$,
 changed  by the production of a substance. The production of the
 substance is given by another function $c_j: \R \rightarrow \R$.

They make the following assumptions:
\begin{itemize}
\item[(1)] The activity of a gene is determined by the state
$b_j\in B_j$ of transcription factor binding sites in its promoter
regions. \item[(2)] Each binding site  can be in one of a
\emph{finite number of states}, characterized by having or not
having bound a particular transcription factor ($B_j$ is a finite
set). \item[(3)] The state $b_{jk}$ of a binding site $B_{jk}$
depends on the concentration $c_{j}(t)$  of the respective
transcription factors. \item[(4)] Depending on the state $b_j$ of
the binding site $B_j$ a gene can either  be silent or have a
particular activity level. \item[(5)] If a gene $G_j$ is active,
the concentration $c_{j}(t)$  of the substance $i_{j}$ that it
produces  is \emph{linearly growing} with a particular rate,
otherwise it is decreasing (or stays at $0)$.
\end{itemize}

Their multistate generalization is the following:

\begin{itemize}
\item[(a)] A binding site can competitively bind more than one
substance and therefore can have more than two states.
\item[(b)] A gene can have more than two levels of activity.
\item[(c)] A control function  is not a boolean function,
but a mapping which maps a vector of integers into an integer.
\end{itemize}

In FSLM  model, multiple transcription factors can act on several
binding sites to produce a finite output state for a gene. The
finite output state translates to a particular growth or decrease
rate (real valued) for a gene product. Transcription factors are
gene products like any other gene product, and are measured with a
real-valued concentration ($c_{j}(t)$). The concentration of
transcription factors determines the finite state of each binding
site. Time in the  model is continuous, but measurements are made
at a finite number of discrete intervals. The measurements that
are performed is of the concentration of the gene products.

  This model is a simplification of the true
biological process in which the RNA produced by transcription is
later translated into proteins which have their own rate of decay.
Proteins and other cellular species can also interact and activate
or deactivate each other besides interacting with the binding
sites.

This  model has two aspects. One is discrete, given by the control
functions. But, the production of the substance at a given rate is
continuous.

\section{ Sequential Dynamical Systems}
For better understanding the next section  we recall some
definitions of  graphs and sequential dynamical systems
\cite{LP2}.

Let $X$ be a set. Let $\wp _2(X)$ be the set of all two-element
subsets of $X$.
\begin{definition} A (loop free,undirected, finite) graph
$G=(V_G, E_G)$ consists of a finite set $V_G$ of vertices  and a
subset $E_G\subseteq \wp _2(X)$ of edges.
\end{definition}

Let $G$ be a graph. A $1$-neighborhood $N(A)$ of a vertex $a\in
V_G$ is the set
\[ N(a):=\{ b\in V_G| \{a,b\}\in E_G\ or \ a=b\}.\]

 Let $V_G=\{a_1,\ldots ,a_n\}$. Let $(k[a_i],a_i\in V_G)$ be a family
of sets.  Define
\[k^n:= k[a_1]\times \ldots k[a_n]=\prod _{a_i\in V_G}k[a_i],\]
the set of global states of $G$.
\begin{definition}\label{local}
A function $f:k^n\rightarrow k^n$ is called local at $a_i\in V_G$
if
\[f(x_1,\ldots ,x_n)=(x_1,\ldots ,x_{i-1}, f^i(x_1, \ldots , x_n),
x_{i+1}, \ldots ,x_n),\] where $f^i(x_1, \ldots , x_n)\in k[a_i]$
depends only on the variables  in the $1$-neighborhood $N(a_i)$ of
the vertex $a_i$.
\end{definition}
\begin{definition}
A sequential dynamical system (SDS), $\F=(Y, (k[a_i]),(f_i),\alpha
)$ consists of
\begin{enumerate}
\item a finite graph $Y$ with $n$ vertices, \item a family of sets
$(k[a_i],a_i\in V_Y)$ in $Z$, \item a family of local functions
$(f_i: k^n \rightarrow k^n,\ where \ f_i \ local \ at \  a_i)$,
\item and a word $\alpha =\alpha _Y=(\alpha _1,\ldots ,\alpha
_r)\in V_Y^*$ in the Kleene closure of the set of vertices $V_Y$,
called an update schedule ( i.e. a map $\alpha :\{1,\ldots
,r\}\rightarrow V_Y)$.
\end{enumerate}
\end{definition}
The world $\alpha$ is used to define the global update function of
an SDS as the function \[ F= f_{\alpha _r} \circ \cdots \circ
f_{\alpha _1}:k^n\rightarrow k^n.\]

The length of the update schedule $\alpha =(\alpha _1, \cdots
,\alpha _r)$ is $r$. The global update function of an SDS defines
its dynamical behavior, properties of limit cycles, transients,
etc..
\begin{definition}
A parallel dynamical system or a finite dynamical system is a
function $F:k^n\rightarrow k^n$.
\end{definition}
\begin{remark}
 Every parallel system can be represented as a sequential system by
doubling the number of nodes and first copying the old states to
the new variables. Conversely, after compose all the local update
functions in a sequential system, then the global update function
$F:k^n\rightarrow k^n$ has coordinate functions (different from
the local update functions in general), and we can think of the
system as being a parallel system given by the coordinate
functions.  So the two representations are equivalent. If the
system one wants to model is naturally sequential, then the
representation as an SDS is generally better because it makes
important system properties explicit.
\end{remark}

\section{Two models}
In this section we define  two models for genetics networks. First
we introduce the definition of a translated function.

\begin{definition}
A vector function $\delta_{fg}=(\delta _1, \ldots , \delta _n ):
\R ^n  \rightarrow \Z ^n $ is a translated function between two
vector functions $g=(g_1,\ldots , g_n):\R ^n \rightarrow \R ^n$,
and $f=(f_1, \ldots , f_n):\Z^n \rightarrow \Z^n $ if:
\[\delta _i \circ g_i= f_i \circ \delta _i , \hbox{ for all } i=1,\ldots , n .\]
That is, if the following diagram commute
\[\begin{array}{lll}
\R ^n & \rightarrow^{g} & \R^n\\
\delta \downarrow & & \delta \downarrow \\
\Z ^n & \rightarrow^{f} & \Z ^n
\end {array}\]
\end{definition}

\begin{definition}  A finite state model (FSM)
$\Gamma = \{Y, (B_j), (F_j), (c_j), \delta _{Fc} \}$, consists of
\begin{enumerate}
\item a finite graph $Y$ with $n$ vertices, $Y$  is the supported
graph of relations between the $n$ genes $G_j$ with vertices
$V_Y=\{g_1,. . .,g_n\}$, \item a family of finite sets $B_j$
(binding sites), for each $g_j\in V_Y$, \item a vector function
$F=(F_1, \ldots , F_n): \prod _{j=1}^n B_j \rightarrow  \prod
_{j=1}^n B_j $ such that $F_j: \prod _{j=1}^n B_j\rightarrow B_j$,
\item a vector function $c=(c_1,\ldots, c_n): \R ^n \rightarrow \R
^n$, where $\R$ is the set of real number, \item a translated
function $\delta _{Fc}$ between the two vector functions.
\end{enumerate}
\end{definition}

\begin{definition}
Let $\{t_0,t_1,\cdots ,t_n\}$ be a set of real numbers, such that
$t_0<t_1<\cdots <t_n$. Let $(a_0,b_0), (a_1,b_1), \ldots ,
(a_n,b_n)$ be $n+1$ pair of real numbers. Suppose that

\[\begin{array}{ccc}
a_0t_1+b_0 & = & a_1t_1+b_1\\
a_1t_2+b_1 & = & a_2t_1+b_2\\
&\cdots &\\
a_{n-2}t_{n-1}+b_{n-2} & = & a_{n-1}t_{n-1}+b_{n-1}
\end{array}\]
We call the function
\[c(t)=\left \{\begin{array}{lcc}
0& if & t<t_0 \ or \ t>t_n \\
a_it+b_i & for & t\in [t_i,t_{i+1}], i=0,1,\ldots , n-1
\end{array} \right .\]
a sectional linear function.
\end{definition}

As a consequence of section $2$ we have proved part of the
following theorem.
\begin{theorem} Brazma-Schlitt Model is a finite state  model, where the
functions $c_j$ are sectional linear functions.
\end{theorem}
\begin{proof}
Here, we only need to see that the functions  $c_j$ give the
concentration  of the substance.
\end{proof}
We assume the following:
\begin{enumerate}
  \item there are  $n$ genes in the network $\textit{N}$,
  \item for each  $1\leq j\leq n$ we have $m_j\in \Z^+$ binding sites
  $\{B_{j1},\ldots ,B_{jm_j}\}$,
  \item $B_{jk}$ is a finite set, and $B_{jk}\subset  \Z$, for all
$j,$
  and $k$,
  \item one gene can be interact with another gene, and we describe
this   situation  by a graph $Y$, with set of vertices
$V_Y=\{g_1,\ldots,g_n\}$,
  \item the environment of the network $\textit{N}$ is the set\[B=\prod
_{j=1}^n
  B_j, \hbox { where } B_j=B_{j1}\times \cdots \times B_{jm_j},\]
  \item for each  gene $g_j$ we have a  local function, $F_j:B\rightarrow B$, in the
sense
  of Definition \ref{local}.
\end{enumerate}

For all function $F:\Z ^n \rightarrow \Z^n$, we define the
coordinated functions $F_i$ as follows: \[ F(x_1, \ldots , x_n
)=(F_1(x_1, \ldots , x_n), \ldots , F_n(x_1, \ldots , x_n)).\]

\begin{definition}\label{gsds}
A \emph{genetic sequential dynamical system} (GSDS) consists of
$\F=(Y, (B_j),(f_j),(c_j),\alpha , \delta)$, where
\begin{enumerate}
\item $Y$ is the support graph of relations between genes with
vertices $V_Y=\{g_1,\ldots ,g_n\}$, \item $B_j$ is a finite set,
for all $j$, and $B=\prod _{j=1}^n B_j$,
 \item a family of local functions $f_j:B\rightarrow B$, (genetic functions),
\item a word $\alpha$ with the order of interaction of functions,
that is a function $F=f_{\alpha(n)} \circ \cdots \circ f_{\alpha
(1)}=(F_1, \ldots , F_n): \Z ^n \rightarrow \Z^n$
 \item a vector  function $c=(c_1,\ldots , c_n): \R ^n \rightarrow \R ^n $,
  \item a translated function $\delta_{Fc}$.
 \end{enumerate}
\end{definition}
\begin{definition}
A genetic network is a pair $\Gamma =(F,c)$. A genetic network is
compatible if there exists a translated function $\delta_{Fc}$.
\end{definition}

\begin{theorem}
The genetic sequential dynamical system is a generalization of the
Brazma-Schlitt  model.
\end{theorem}
\begin{proof}
In order to prove the theorem, we see how all the considerations
of Brazma-Schlitt (BS) model are included in the definition of
GSDS.

 The interaction between genes is given by a  graph $Y$ with $n$
vertices
 $g_1,\ \ldots \ ,\ g_n$ and an edge $\{g_i,g_k\}$  if the gene
 $i$ has any relation
with  the gene $k$. So, we have the first condition of Definition
\ref{gsds}.

We can observe that  in  Definition \ref{gen} a gene in action
\emph{is a complex molecular machinery}. So a gene is a function
that \emph{can read information and make a particular type of
protein}.
 Where does a gene
read the information? In the binding sites $B_j$ and the protein
again changes the environment $B_j$. The Brazma model considers
functions $\overline F_j$ from $B_j$ to $B_j$. Since the protein
can make changes in the $1$-neighborhood,  the genetic function
$F_j$ is from $B$ to $B$, and it  is a local function. In the BS
model, the binding sites are finite sets. So, we have conditions
2,  and 3.

The genes act in an order, which implies an order in the
composition of those genetic functions. Thus condition 4 holds.

The inclusion of a family of functions $(c_j)$ and the translated
functions in the definition of GSDS gives the possibility  to see
the continuous and discrete sides of the genetic networks.

Then our claim holds.

\end{proof}
\section{Examples}
In all the examples, we obtained the functions using partially
defined functions and polynomial representation, see \cite{A1}.
\begin{example}
We describe the Boolean model presented in \cite{ITK,O} using a
GSDS. In this example the data is given with discrete values: 0
and 1.
\begin{enumerate}
\item The digraph: $Y \ \  \ \ \ \begin{array}{ccc}
g_1&\rightarrow & g_3\\
& \searrow  & \uparrow \\
g_0 & \rightarrow & g_2
\end{array}$

\item $X=\Z_2=\{0,1\}$.

\item  The functions are the following:  $\begin{array}{ccl}
f_0(x_0,x_1,x_2,x_3)&= &(1,x_1,x_2,x_3) \\
f_1(x_0,x_1,x_2,x_3)&= &(x_0,1,x_2,x_3)  \\
f_2(x_0,x_1,x_2,x_3)&= &(x_0,x_1,x_0x_1,x_3) \\
f_3(x_0,x_1,x_2,x_3)&= &(x_0,x_1,x_2,x_1(x_2+1))
\end{array}$
\item The schedule   $\alpha=\left(
\begin{array}{llll}
0&1& 2&3 \\
3&2& 1& 0
\end{array} \right)$,
\item The global function $f=f_0 \circ f_1 \circ f_2 \circ
f_3:X^4\rightarrow X^4$,
\[f(x_0,x_1,x_2,x_3)= (1,1,x_0x_1,x_1(x_2+1)).\]
\end{enumerate}
We can observe that if we change the order we can not obtain the
same function.
\end{example}

\begin{remark}
If we have only three states for genes we have the finite field
$\Z_3=\{-1,0,1\}$, that is the integers modulo $3$, with $1+1=-1$.
If we have four possible states for genes we can use a finite
field with $4$ elements. A finite field GF$(4)$ can be represented
as: $GF(4)=\{0,1,\alpha, \alpha ^2\}$, where $\alpha$ is a root of
the polynomial $z^2+z+1$, that is $\alpha ^2=\alpha +1$  (with
coefficients in $\Z_2=\{0,1\}$). We denote $0=00,\ 1=01,\
\alpha=10,\ \alpha ^2=11$, then the operations $+$ and $\times$
are the follows:

\[ \begin{array}{ccccc}
+ & 00 & 01& 10 & 11\\
00 & 00 & 01& 10& 11\\
01 & 01 & 10 & 11& 00\\
10 & 10 & 11 & 00 & 01\\
11 & 11 & 10 & 01 & 00
\end{array} \ \ \ \ \ \ \
\begin{array}{ccccc}
\times & 00 & 01& 10 & 11\\
00 & 00 & 00& 00& 00\\
01 & 00 & 01 & 01& 11\\
10 & 00 & 10 & 11 & 01\\
11 & 00 & 11 & 01 & 10
\end{array} \]

\end{remark}
\begin{example}
We describe the  example which appear in \cite{J} of the
generalized logical method developed by Thomas and colleagues
\cite{ T}. Here we use the FSM and in this case it is not linear
and we have the data with two or three values in the integers. We
have three genes, and the regulatory network is the following:
\begin{enumerate}
\item The digraph: $Y \ \  \ \ \ \begin{array}{ccc}
& & g_3 \circlearrowleft\\
 & \nearrow & \downarrow \\
g_1 & \rightleftarrows & g_2
\end{array}$

\item For genes $g_1$ and  $g_3$ we have $\Z_3=\{0,1,2\}$, and for
gene $g_2$ we have $X_2=\{0,1\}$,

\item   The functions are:
\[\begin{array}{ccl}
f_1(x_1,x_2,x_3)&= &-x_2  \\
f_2(x_1,x_2,x_3)&= &1+x_1^2x_3^2 \\
f_3(x_1,x_2,x_3)&=
&2+x_1+2x_3+x_1x_3+2x_1^2+x_3^2+2x_1^2x_3+2x_1x_3^2+x_1^2x_3^2
\end{array}\]
\item The global function \[F= (f_1, f_2, f_3):\Z_3\times
X_2\times \Z_3 \rightarrow \Z_3\times X_2\times \Z_3\] such that
\[F(x_1,x_2,x_3)= (-x_2,1+x_1^2x_3^2,2+x_1+2x_3+x_1x_3+2x_1^2+x_3^2+2x_1^2x_3+2x_1x_3^2+x_1^2x_3^2).\]
\end{enumerate}

\end{example}
\begin{example}
 We assume that by microarray experiment we have the following data,
 and denote the concentration of each gene $j$ by $c_j(t)$:

 \[ \begin{array}{cccc}
 t & G_1 & G_2& G_3\\
 0& c _1(0)=0.5 & c _2(0)=1.2& c_3(0)=0.5\\
 1 &c _1(1)=0.78 & c _2(1)=1.2 & c_3(1)=1.25\\
  2 &c _1(2)=1.5 & c _2(2)=1.2 & c_3(2)=1.5\\
   3 &c _1(3)=0.5 & c _2(3)=1.2 & c_3(3)=0.5
  \end{array}\]
  The vector of concentrations  is
  $c (t)=(c_1(t),c_2(t),c_3(t)),$
 and if we suppose that the functions $c_j$ are sectional linear
 functions then
\[c (t)= \left\{ \begin{array}{ll}
(0.28t+0.5,1.2,0.75t+0.5),& \ when \ t\in
[0,1]\subset \R \\
(0.72t+0.06,1.2,0.25t+1), & \ when \ t\in [1,2]\subset
\R \\
(-t+3.5,1.2,-t+3.5), & \ when \ t\in [2,\infty)\subset \R
\end{array} \right\}\]
On the other hand, we will assume that the average concentration
level of $G_1$ is $0.78$,  $G_2$ is
 $0.75$, and $G_3$ is $1.5$. So, we give values to the states of genes $G_1$, $G_2$, and
$G_3$.
\[\begin{array}{cccc}
\delta_1: & \{\ less \ than \  0.78 \mapsto -1, & 0.78\mapsto 0,& \ more \ than \ 0.78 \mapsto 1\}\\
\delta_2: &\{ \ less \ than \  0.75 \mapsto -1,& 0.75\mapsto 0 & \ more \ than \ 0.75 \mapsto 1 \}\\
\delta_3: & \{ \ less \ than \  1.25 \mapsto -1,& 1.25\mapsto 0, &
\ more \ than \  1.25 \mapsto 1\}
\end{array}\]
We suppose that we have for each gene $G_j$ a binding site
$B_j=X=\{-1,0,1\}$,  and we consider  the operations in the finite
field $X=\Z_3$. Our problem is the following: we know $c (t)$ by
microarray experiment, we suppose the vector $c (t)$ is a vector
of sectional affine functions, but we want to obtain a function
$f$ such that $f(-1,1,-1)=(0,1,0)$ for $t=1$, $f(0,1,0)=(1,1,1)$
for $t=2$, and $f(1,1,1)=(-1,1,-1)$ for $t=3$. In this case, one
of the possible functions is
$f(x_1,x_2,x_3)=(x_1+x_2,x_2,x_3+x_2)$. In addition, we can obtain
a graph $Y$ if we observe how the genes change with the action of
$F$.
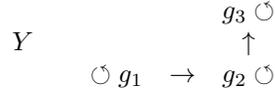
\[Y \ \  \ \ \ \begin{array}{ccc}
& & g_3 \circlearrowleft \\
& & \uparrow \\
\circlearrowleft g_1 & \rightarrow & g_2 \circlearrowleft
\end{array}\]
Now, we have a GSDS $\F=(Y, \Z_3 , \{f_j\} , (\delta _j),
(c_j),\alpha )$ for this dataset:
\begin{enumerate}
\item a collection $x_1,x_2,x_3$ of variables, which  take on
values in a finite field $X=\Z_3$

 \item  $Y$ is the support directed graph of relations between genes
with vertices $\{g_1,g_2,g_3\}$,

 \item for each $j=1,2,3, $ the local update
 functions \[f_1(x_1,x_2,x_3)=(x_1+x_2,x_2,x_3), \
 f_2(x_1,x_2,x_3)=(x_1,x_2,x_3), \]
  \[ \hbox {  and } f_3(x_1,x_2,x_3)=(x_1,x_2,x_2+x_3).\]

\item a  schedule  $\alpha =\left(
\begin{array}{lll}
1& 2&3 \\
1& 2& 3
\end{array} \right)$,

\item the global function $f=f_3\circ f_2 \circ f_1:X^3\rightarrow
X^3$ obtained by the schedule $\alpha $,

 \item a vector function $c(t)=(c_1(t),c_2(t),c_3(t))$,
\item the translated functions $\delta_i$.
\end{enumerate}
\end{example}

\end{document}